\documentclass[a4paper,12pt]{amsart}
\usepackage{fancyheadings,amsmath,amscd,amsthm,amsfonts,latexsym,mathrsfs,amssymb}
\DeclareMathOperator{\Div}{div}
\DeclareMathOperator{\trace}{trace}
\DeclareMathOperator{\grad}{grad}

\DeclareMathOperator{\RicM}{^{\it M}\!Ricci}
\DeclareMathOperator{\RicN}{^{\it N}\!Ricci}
\DeclareMathOperator{\RicP}{^{\it P}\!Ricci}
\DeclareMathOperator{\dif}{d}

\DeclareMathOperator{\Lie}{\mathcal{L}}

\renewcommand{\H}{\mathcal{H}}
\newcommand{\V}{\mathcal{V}}
\newcommand{\F}{\mathcal{F}}
\DeclareMathOperator{\Bh}{{\it B}^{\H}}
\DeclareMathOperator{\Bv}{{\it B}^{\V}}
\DeclareMathOperator{\M}{\overset{\it M}{\nabla}}
\DeclareMathOperator{\N}{\overset{\it N}{\nabla}}

\def \b{\beta}

\def \G{\Gamma}
\def \l{\lambda}
\def \O{\Omega}
\def \phi{\varphi}
\def \s{\sigma}
\def \D{\Delta}

\def \Co{\mathbb{C}\,}

\begin{document}

\title{Harmonic morphisms with one-dim\-en\-sion\-al fibres\\
on Einstein manifolds}
\author{Radu Pantilie and John C.\ Wood} 
\thanks{The authors gratefully acknowledge that this work was 
done under E.P.S.R.C. grant number GR/N27897} 
\email{r.pantilie@leeds.ac.uk\,, j.c.wood@leeds.ac.uk\,.}

\newtheorem{thm}{Theorem}[section]
\newtheorem{lem}[thm]{Lemma}
\newtheorem{cor}[thm]{Corollary}
\newtheorem{prop}[thm]{Proposition}

\theoremstyle{definition}
\newtheorem{defn}[thm]{Definition}
\newtheorem{rem}[thm]{Remark}
\newtheorem{exm}[thm]{Example}

\numberwithin{equation}{section}

\maketitle
\thispagestyle{empty} 
\vspace{-0.5cm}
\section*{Abstract}
\begin{quote}
{\footnotesize
We prove that, from an Einstein manifold of dimension greater than or 
equal to five, there are just two types of harmonic morphism with 
one-dim\-en\-sion\-al fibres. This generalizes a result of R.L.~Bryant 
who obtained the same conclusion under the assumption that the 
domain has constant curvature.}
\end{quote}

\section*{Introduction} 

\indent
Harmonic morphisms between Riemannian manifolds are smooth maps which preserve Laplace's 
equation. By a basic result of B.~Fuglede \cite{Fug} and T.~Ishihara \cite{Ish}\,, they 
can be characterised as harmonic maps which are horizontally weakly conformal. See \cite{BaiWoo2} 
for a general account, and \cite{Gudbib} for a frequently updated bibliography and other useful  
information on harmonic morphisms.\\ 
\indent
The first classification results for harmonic morphisms with one-dim\-en\-sion\-al fibres 
were obtained by P.~Baird and J.C.~Wood (see \cite{BaiWoo2}). In \cite{Bry}\,, R.L.~Bryant 
proved that, from a constant curvature Riemannian manifold of dimension greater than or 
equal to four, there exist only two types of harmonic morphisms with one-dim\-en\-sion\-al fibres (see \cite{BaiWoo2} and \cite{Pan-thesis} for alternative proofs). The two types can be nicely 
described in terms of the geometrical properties of the foliations formed by the components 
of the regular fibres. For this, recall \cite{Woo}\,, \cite{Pan-thesis} that we say 
that a foliation \emph{produces harmonic morphisms} if it can be locally defined by 
submersive harmonic morphisms. Then the two types of foliations which produce harmonic 
morphisms which appear in Bryant's theorem are the following (note that in \cite{Bry} the 
ordering is different):\\  
\indent
{\bf Type 1}. Riemannian one-dim\-en\-sion\-al foliations locally generated by 
Killing vector fields \cite{Bry}\,.\\ 
\indent
{\bf Type 2}. Homothetic foliations by geodesics orthogonal to an umbilical 
foliation by hypersuraces \cite{BaiEel}\,.\\  
Accordingly, a harmonic morphism with one-dim\-en\-sion\-al fibres is said to be of type 1 (respectively, 
type 2) if the components of its regular fibres form a foliation of type 1 (respectively, type 2). 
Note that the gradient of the dilation at regular points is horizontal for type 1 and vertical for 
type 2. Also, recall \cite{Bai} that a harmonic morphism with one-dim\-en\-sion\-al fibres is always 
submersive if the domain has dimension at least five, whilst if the domain is of dimension four,  
then the set of critical points is discrete (see \cite{BaiWoo2} for the case when the domain has dimension  three).\\   
\indent
We prove (Theorem \ref{thm:dmain}) that any harmonic morphism with one-dim\-en\-sion\-al fibres from an 
Einstein manifold of dimension greater than or equal to five is either of type 1 or of type 2, thus 
generalizing the result of \cite{Bry}\,. Note that, from an 
Einstein \emph{four}-manifold, there is just one more type of harmonic morphism with 
one-dim\-en\-sion\-al fibres \cite{Pan-thesis}\,, \cite{Pan-4to3}\,.\\ 
\indent
Section 1 contains some generalities on foliations which produce harmonic morphisms. 
The main classification result is obtained in Section 2 (Theorem \ref{thm:dmain}). 
If $(M^{n+1},g)$ is an Einstein manifold and $\phi:(M^{n+1},g)\to(N^n,h)$ is a harmonic 
morphism of type 2 then $(N^n,h)$ is also an Einstein manifold. (This follows, for example, from  well-known results on warped products (see \cite[Chapter 9]{Bes}\,).) However, this is not 
true for a harmonic morphism of type 1. In Section 3 we refine the classification of harmonic morphisms with one-dim\-en\-sion\-al fibres in the case when \emph{both} the domain and codomain are  
Einstein manifolds. The result (Theorem \ref{thm:dbetween}) shows that only two 
well-known constructions (see \cite[Chapter 9]{Bes}\,), illustrated by the Hopf fibrations and 
suitable warped products over Einstein manifolds, can occur. 
In the Appendix we prove a fact on horizontally holomorphic submersions 
which is needed in the proof of Theorem \ref{thm:dbetween}\,.

\section{Foliations which produce harmonic morphisms}

\indent
Let $\V$ be (the tangent bundle of) a smooth foliation on a smooth Riemannian manifold 
$(M,g)$ and set $\H=\V^{\perp}$\,. Recall (see \cite{Pan-thesis}) that we say that $\V$ 
\emph{produces  harmonic morphisms} if it can be locally defined by submersive harmonic morphisms. 
A foliation of codimension two produces harmonic morphisms if and only if it is a conformal foliation 
with minimal leaves \cite{Woo}\,.  
For foliations of codimension not equal to two which produce harmonic morphisms we 
have the following characterisation of R.L.~Bryant \cite{Bry} (see also \cite{BaiWoo2} and 
\cite{Pan-thesis}\,). 
\begin{prop} \label{prop:Bry} 
Let $\V$ be a foliation of ${\rm codim}\,\V=n\neq2$ on the Riemannian manifold  $(M,g)$\,. Then 
the following assertions are equivalent:\\ 
\indent
{\rm (i)} $\V$ produces harmonic morphisms;\\ 
\indent
{\rm (ii)} $\V$ is a conformal foliation and the one-form $$\b=(n-2)\trace(\Bh)^{\flat}-n\trace(\Bv)^{\flat}$$ 
is closed where $\H=\V^{\perp}$, and $\Bh$ (respectively, $\Bv$) denotes the second fundamental 
form of $\H$ (respectively, $\V$). (As is usual, we denote by the same letter $\H$ 
(respectively, $\V$) 
the horizontal (respectively, vertical) distribution and the orthogonal projection onto it.) 
\end{prop} 
\begin{proof}
If $\V$ is a conformal foliation then it is easy to see that, locally, we have 
\begin{equation*} 
\trace(\Bh)=n\,\V(\grad\log\l) 
\end{equation*} 
for any local dilation $\l$ of $\V$ (see, for example, \cite[Chapter 2]{BaiWoo2}).\\ 
\indent
{}From a fundamental formula of P.~Baird and J.~Eells \cite{BaiEel}\,, it follows that 
a conformal foliation $\V$ produces harmonic morphisms if and only if, in the 
neighbourhood of each point, there 
exists a local dilation $\l$ such that
\begin{equation} \label{e:BairdEells}  
\trace(\Bv)=-(n-2)\H(\grad\log\l)\;. 
\end{equation}
Therefore, for any such local dilation $\l$\,, we have $\b=n(n-2)\dif\log\l$ and the proof follows.  
\end{proof} 

\indent
We shall also need the following definition \cite{Pan-thesis}. 

\begin{defn} \label{defn:homfoln} 
We say that a foliation is \emph{homothetic} if it can be locally defined by horizontally 
homothetic submersions (i.e., horizontally conformal submersions whose dilations are constant 
along horizontally curves). 
\end{defn} 

\begin{rem} \label{rem:homfoln} 
1) Obviously any Riemannian foliation is homothetic. Also, a homothetic foliation is conformal.  Conversely, a conformal foliation is homothetic if and only if the mean curvature form 
$\trace(\Bh)^{\flat}$ of its orthogonal complement 
$\H$ is closed \cite{Pan-thesis}\,. \\ 
\indent 
2) If $\V$ is a Riemannian foliation of codimension $n$ then 
$$\b=-\,n\trace(\Bv)^{\flat}\;.$$ 
By Proposition \ref{prop:Bry}\,, if $n\neq2$\,, $\V$ produces harmonic morphisms if and only 
if it has closed mean curvature form. Further, if $\V$ is a \emph{one-dim\-en\-sion\-al} Riemannian 
foliation, then this is easily seen to be equivalent to the fact that $\V$ is locally generated by 
Killing vector fields. Therefore \emph{a Riemannian one-dim\-en\-sion\-al foliation of 
codimension not equal to $2$ 
produces harmonic morphisms if and only if it is locally generated by Killing vector fields} \cite{Bry} (see also \cite{BaiWoo2}\,, \cite{Pan} and \cite{Pan-thesis}\,).\\ 
\indent
3) If $\V$ is a foliation with minimal leaves of codimension $n$ then $$\b=(n-2)\trace(\Bh)^{\flat}\;.$$  
{}From Proposition \ref{prop:Bry} it follows that \emph{a foliation with minimal leaves of 
codimension not equal to $2$ produces harmonic morphisms if and only if it is a homothetic 
foliation} (this is essentially a result of \cite{BaiEel}\,). 
\end{rem}

\indent 
The following result answers a question of P.~Baird. 

\begin{prop} \label{prop:analytic}
Let $(M,g)$ be a real-analytic Riemannian manifold, and let $\V$ be a foliation which 
produces harmonic morphisms on $(M,g)$.\\ 
\indent
Then $\V$ is a real-analytic foliation. Moreover, if ${\rm codim}\,\V\neq2$\,, then any harmonic 
morphism produced by $\V$ is a real-analytic map onto a real-analytic Riemannian manifold . If  
${\rm codim}\,\V=2$\,, then $\V$ is locally defined by real-analytic submersive harmonic morphisms 
onto real-analytic Riemannian two-manifolds. 
\end{prop}  
\begin{proof} 
Suppose that ${\rm codim}\,\V=n$\,, and let $\phi:(U,g|_U)\to(N,h)$ be a submersive 
harmonic morphism whose fibres are open subsets of leaves of $\V$\,. Let $\xi:(V,h|_V)\to\mathbb{R}^n$ 
be local coordinates on $(N,h)$ given by harmonic functions 
(see \cite{Ish} or \cite{BaiWoo2} for a proof that such local coordinates exist). Then, because $\phi$ is a harmonic morphism, 
$\xi\circ\phi:(\phi^{-1}(V),g|_{\phi^{-1}(V)})\to\mathbb{R}^n$ is a submersive harmonic 
map. {}From the fact that $(M,g)$ is real-analytic it follows that the components of 
$\xi\circ\phi$ are real-analytic functions on $M$\,, and hence $\V$ 
is locally defined by real-analytic submersions. Thus $\V$ is a real-analytic 
foliation.\\ 
\indent
If $n=2$\,, the proof is ended.\\  
\indent
Suppose that $n\neq2$\,. {}From the fact that both $(M,g)$ and $\V$ are analytic it follows  
that the one-form $\b$ is analytic and hence 
the dilation of any harmonic morphism produced by $\V$ is an analytic function.\\ 
\indent
Let $\phi:(U,g|_U)\to(N,h)$ be a submersive harmonic morphism produced by $\V$\,. We 
shall show that $(N,h)$ is real-analytic in harmonic coordinates. Indeed, since the dilation 
$\l$ of $\phi$ is analytic and $\phi^*(h)=\l^2g|_{\H}$ we have that, in harmonic coordinates,  $h$ has analytic components. Then, the analyticity of the atlas 
$\mathcal{A}$ on $N$ given by the harmonic coordinates follows from the analyticity 
of harmonic functions on analytic Riemannian manifolds (applied to the charts 
of $\mathcal{A}$). 
\end{proof}  

\begin{cor} 
Let $\phi:(M,g)\to(N,h)$ be a submersive harmonic morphism from a real-analytic manifold onto a smooth manifold.\\ 
\indent
If $\dim N\geq3$ then there exists a real-analytic structure on $N$ with respect to 
which $\phi$ and  
$h$ are real-analytic. If $\dim N=2$ then $\phi$ is real-analytic with respect to the real-analytic structure on $N$ induced by the conformal structure of $h$\,.  \qed
\end{cor}

\section{The proof of the main theorem}

\indent
We shall need two straightforward lemmas. For the reader's convenience we include their 
proofs here. 

\begin{lem} \label{lem:1} 
For $n\geq1$, let $P:\mathbb{C}^n\times\mathbb{C}\to\mathbb{C}$ be defined by 
$$P({\bf a},\l)=\l^n+a_1\l^{n-1}+\cdots+a_n\;,\qquad({\bf a}=(a_1,\ldots,a_n)\in\Co^n\,,\: \l\in\mathbb{C}\,).$$ 
Let $\bigl({\bf a}^{(k)}\bigr)\subseteq\mathbb{C}^n$ be a convergent sequence. \\ 
\indent
Then the set 
$S=\bigl\{\,\l\in\mathbb{C}\;|\;\exists\,k\,:\,P({\bf a}^{(k)},\l)=0\,\bigr\}$ is bounded. 
\end{lem}
\begin{proof} 
Let $\widetilde{P}:\mathbb{C}^n\times\mathbb{C}P^1\to\mathbb{C}P^1$ be defined by 
$$\widetilde{P}({\bf a},[\l_0,\l_1])=[\l_0^n,\l_1^n+a_1\l_0\l_1^{n-1}+\cdots+a_n\l_0^n]\;,\qquad 
({\bf a}\in\mathbb{C}^n\,,\:[\l_0,\l_1]\in\mathbb{C}P^1\,).$$  
\indent
Obviously 
$S=\bigcup_k S_k$ where 
$S_k=\bigl\{\,\l\in\mathbb{C}\;|\;\widetilde{P}({\bf a}^{(k)},[1,\l]\,)=[1,0]\,\bigr\}$\,;   
suppose that $S$ is unbounded. For each $k$\,, let $\l^{(k)}\in S_k$ 
be such that 
$|\l|\leq|\l^{(k)}|$ for every $\l\in S_k$. Because $S$ is unbounded we must have that 
$\bigl(\l^{(k)}\bigr)$ is unbounded. By passing to a subsequence, if necessary, we can suppose 
that $\l^{(k)}\to\infty$\,.\\ 
\indent
Because $\l^{(k)}\in S_k$ we have $\widetilde{P}({\bf a}^{(k)},[1,\l^{(k)}]\,)=[1,0]$. But 
$$\lim_{k\to\infty}\widetilde{P}({\bf a}^{(k)},[1,\l^{(k)}]\,)=\widetilde{P}({\bf b},[0,1])=[0,1]$$  
where ${\bf b}=\lim_{k\to\infty}{\bf a}^{(k)}$\,, a contradiction. Hence $S$ must be bounded.  
\end{proof}

\indent
For the second lemma we need the following definition. 

\begin{defn} \label{defn:consdiag} 
Let $E\to M$ be a vector bundle endowed with a Riemannian 
metric $h$\,, and let $T\in\G(\odot^2E^*)$\,. Say that $T$ can be \emph{consistently diagonalized} 
at $x_0\in M$ if there exists 
an open neighbourhood $U$ of $x_0$ and a local orthonormal frame $\{e^j\}$ 
on $U$ for $E^*$ such that $$T=\sum_{j,k}\mu_j\,\delta_{jk}\,e^j\otimes e^k$$  
for some smooth functions $\mu_j:U\to\mathbb{R}$\,.\\ 
\indent
A similar definition can be made for a field of self-adjoint endomorphisms 
$\widetilde{T}\in\G({\rm End}\,E)$\,. 
\end{defn} 

\begin{lem} \label{lem:consdiag} 
Let $E\to M$ be a vector bundle endowed with a Riemannian 
metric $h$\,, and let $T\in\G(\odot^2E^*)$\,.\\ 
\indent
Then $T$ can be consistently diagonalized at each point of a dense open subset of $M$\,. 
\end{lem} 
\begin{proof} 
The proof is by induction on the rank (fibre dimension) of $E$.\\ 
\indent
For ${\rm rank}\,E=1$ the lemma is trivial.\\ 
\indent
Suppose that the assertion of the lemma is true for ${\rm rank}\,E<n$\,; we shall prove that 
the assertion is true for ${\rm rank}\,E=n$\,.\\ 
\indent
Let $P(x,\l)=P_x(\l)$ be the characteristic polynomial of $T_x$\,, 
with respect to $h_x$ ($x\in M$). For $p=1,\ldots,n$\,, set 
$$G_p=\bigl\{\,x\in M\;|\;P_x\:{\rm has\:a\:root\:of\:order\:at\:most\:}p\,\bigr\}\,,$$ 
and set $G_0=\emptyset$\,. Because 
\begin{equation*} 
M=\overset{n}{\underset{p=1}{\bigcup}}\,\overline{G}_p\setminus\overline{G}_{p-1}
\subseteq\overset{n}{\underset{p=1}{\bigcup}}\,\overline{G_p\setminus\overline{G}_{p-1}} 
=\overline{\overset{n}{\underset{p=1}{\bigcup}}\,G_p\setminus\overline{G}_{p-1}}\subseteq M\,,  
\end{equation*}
we have that 
\begin{equation*} 
M=\overline{\overset{n}{\underset{p=1}{\bigcup}}\,G_p\setminus\overline{G}_{p-1}} 
\end{equation*} 
where $\overline{A}$ denotes the closure of the set $A$\,. To complete 
the proof, it suffices to prove that each 
$x\in\overset{n}{\underset{p=1}{\bigcup}}\,G_p\setminus\overline{G}_{p-1}$ 
has an open neighbourhood $U$ such that $E|_U=E_1\oplus E_2$ where $E_1$ and $E_2$ are complementary orthogonal vector subbundles of $E$ 
of positive rank such that $T|_{E_1\otimes E_2}=0$\,.\\ 
\indent
Let $p\in\{1,\ldots,n\}$ be such that $G_p\setminus\overline{G}_{p-1}\neq\emptyset$ and 
let $x_0\in G_p\setminus\overline{G}_{p-1}$\,. Let $\l_0$ be a root 
of $P_{x_0}$ of order at most $p$\,. Because $x_0$ is not in $\overline{G}_{p-1}$ we have that 
$\l_0$ has order $p$\,. Then, by the smooth version of the Weierstrass Preparation Lemma, in an 
open neighbourhood $U\subseteq M\setminus\overline{G}_{p-1}$ of $x_0$  we have 
$P(x,\l)=Q(x,\l)\,R(x,\l)$ where $Q$ is a polynomial of degree $p$ in $\l$ such that 
$Q(x_0,\l)=(\l-\l_0)^p$ and $R(x_0,\l_0)\neq0$. {}From the fact that $P$ and $Q$ are both 
polynomials in $\l$ (with coefficients smooth functions of $x$\,), it follows that $R$ is also   polynomial in $\l$\,.\\ 
\indent 
We shall show that there exists an open neighbourhood $V\subseteq U$ of $x_0$ such that,  
for each $x\in V$\,, $Q_x$ has a root of order $p$\,. Suppose not. Let 
$\bigl(x^{(k)}\bigr)\subseteq U$ be such that 
$\lim_{k\to\infty}x^{(k)}=x_0$ and, for each $k$\,,  
there exists $\mu^{(k)}\in\mathbb{R}$ such that $Q\bigl(x^{(k)},\mu^{(k)}\bigr)=0$ with 
$\mu^{(k)}$ a root of $Q_{x^{(k)}}$ of order less than $p$\,. Obviously, $\mu^{(k)}$ is also 
a root of $P_{x^{(k)}}$ and, because $x^{(k)}\in U\subseteq M\setminus\overline{G}_{p-1}$\,, 
$\mu^{(k)}$ is a root of order at least $p$ of $P_{x^{(k)}}$\,. Hence $R\bigl(x^{(k)},\mu^{(k)}\bigr)=0$\,. Now, by Lemma \ref{lem:1}, the sequence  $\bigl(\mu^{(k)}\bigr)$ 
is bounded and hence, by passing to a subsequence, if necessary, 
we can suppose that $\lim_{k\to\infty}\mu^{(k)}=\mu_0$ with $\mu_0\in\mathbb{R}$\,. 
Then $R(x_0,\mu_0)=\lim_{k\to\infty}R\bigl(x^{(k)},\mu^{(k)}\bigr)=0$ and also 
$Q(x_0,\mu_0)=\lim_{k\to\infty}Q\bigl(x^{(k)},\mu^{(k)}\bigr)=0$\,. Because $R(x_0,\l_0)\neq0$ 
we have that $\mu_0\neq\l_0$\,. But this implies that $\l_0$ is not a root of order $p$ of 
$Q_{x_0}$\,. It follows that, in an open neighbourhood $V$ of $x_0$\,, we have that $Q_x$ has 
only roots or order $p$ for any $x\in V$. Thus we can write $Q(x,\l)=\bigl(\l-\mu(x)\bigr)^p$\,,   
($(x,\l)\in V\times\mathbb{R}$), where $\mu(x)$ is the root of 
$\partial^{p-1}Q/\partial\l^{p-1}(x,\cdot)$ so that $\mu$ is smooth. Hence 
$$P(x,\l)=\bigl(\l-\mu(x)\bigr)^p\,R(x,\l)\;,\qquad((x,\l)\in V\times\mathbb{R}\,).$$ 
Moreover, because $\partial^{p}P/\partial\l^{p}(x_0,\l_0)\neq0$ we can suppose that 
$\partial^{p}P/\partial\l^{p}\bigl(x,\mu(x)\bigr)$ is non-zero for any $x\in V$\,. 
It follows that $\mu(x)$ is an eigenvalue of order $p$ for $T_x$ for any $x\in V$\,.\\ 
\indent
Let $(E_1)_x$ be the eigenspace of $\mu(x)$ and let $(E_2)_x$ be its orthogonal complement. 
It is easy to see that $E_j=\underset{x\in V}{\bigcup}(E_j)_x$\,, $j=1,2$\,, are smooth 
subbundles of $E$ which have the required properties. The lemma follows. 
\end{proof} 

\indent
Let $\phi:(M^{n+1},g)\to(N^n,h)$ be a submersive harmonic morphism between Riemannian 
manifolds of dimension $n+1$ and $n$\,, respectively. 
Let $\V$ be the foliation formed by the components of the fibres of $\phi$\,, and let $V\in\G(\V)$ 
be a local vertical vector field such that $g(V,V)=\l^{2n-4}$\,. 
Then let $\theta$ 
be the vertical dual of $V$ (i.e., $\theta(V)=1$ and ${\rm ker}\,\theta=\H$ with $\H=\V^{\perp}$). 
Note that, at each point, both $V$ and $\theta$ are uniquely determined up to sign. 
Therefore, we can globally define a Riemannian metric $\overline{h}$ on $M^{n+1}$ by 
$$\overline{h}=\phi^*(h)+\theta^2\;.$$ 
Note that $\phi:(M,\overline{h})\to(N,h)$ is a Riemannian submersion with geodesic fibres. It follows 
that $V$ is an infinitesimal automorphism of $\H$\,; equivalently,  
\begin{equation} \label{e:Vinfautom} 
[V,X]=0\quad{\rm for\: any\: basic\: vector\: field\:}X\in\G(\H)\;. 
\end{equation} 
Set $\O=\dif\!\theta$\,; from \eqref{e:Vinfautom} it follows easily that 
$\O$ is basic, equivalently, $i_V\O=0$ and $\Lie_V\!\O=0$\,. Also,
$\H$ is integrable if and only if $\O=0$\,.\\ 
\indent
To simplify the notation, from now on we shall denote by the same letter $h$ the induced metric 
$\overline{h}$ on $M$\,. Obviously, $h|_{\H}$ is basic.\\ 

\indent
{}From \cite[Lemma B.1.5]{Pan-thesis} we recall the following lemma. 

\begin{lem} \label{lem:Ricci}
Let $n\geq1$ and let $\phi:(M^{n+1},g)\to(N^n,h)$ be a submersive harmonic morphism between 
Riemannian manifolds of dimension $n+1$ and $n$\,, respectively. 
Let $\l=e^{\s}$ be the dilation of $\phi$\,.\\ 
\indent 
If $\RicM$ denotes the Ricci tensor of $(M,g)$ and $\RicN$ denotes the Ricci 
tensor of $(N,h)$\,, then,
\begin{equation} \label{e:riccixy} \begin{split}
\RicM(X,Y)=(\RicN)(\phi_*X,\phi_*Y)&-\tfrac{1}{2}\,e^{(2n-2)\s}\;h(i_{X}\O,i_{Y}\O)\\
-e^{-2\s}\;(\D{}^{M}\s)\,h(X,Y)&-(n-1)(n-2)X(\s)Y(\s)\:,
\end{split} \end{equation}
\begin{equation} \label{e:riccixv} \begin{split}
\RicM(X,V)=\tfrac{1}{2}\,e^{(2n-2)\s}\,(  ^{h}\dif^{*}\!\O)(X)&+(n-1)e^{(2n-2)\s}\;\O(X,\grad_{h}\s)\\
+\;(n-1)X(V(\s))&-(n-1)(n-2)X(\s)V(\s)\:,
\end{split} \end{equation} 
\begin{equation} \label{e:riccivv} \begin{split}
\RicM(V,V)=(n-2)e^{(2n-4)\s}\D^{M}\s&+2(n-1)V(V(\s))\\
-(3n-4)(n-1)V(\s)^{2}&+\tfrac{1}{4}e^{(4n-4)\s}|\O|_h^{2}\:. 
\end{split} \end{equation}
where $^{h}\!\dif^{*}$ denotes the codifferential on $(M,h)$\,.
\end{lem}  

\indent
We shall also need the following result of \cite{Pan-thesis}\,; for the reader's 
convenience we include a proof. 

\begin{prop} \label{prop:horintorhomoth} 
Let $(M,g)$ be an Einstein manifold of dimension at least $4$\,, and let $\V$ be a  one-dim\-en\-sion\-al 
foliation which produces harmonic morphisms on $(M,g)$\,. Suppose that, 
either, the orthogonal complement $\H$ of $\V$ is integrable, or, $\V$ is a homothetic foliation.\\ 
\indent
Then \emph{either},\\  
\indent 
$(\rm{i})$ $\V$ is a Riemannian foliation locally generated by Killing vector fields,  \emph{or}\\ 
\indent
$(\rm{ii})$ $\V$ is a homothetic foliation by geodesics orthogonal to an umbilical foliation by 
hypersurfaces. 
\end{prop}
\begin{proof} 
Since any Einstein manifold can be given a real-analytic structure (see 
\cite[Theorem 5.26]{Bes}\,), from Proposition \ref{prop:analytic} it follows that it is sufficient to prove that one of the cases (i) or (ii) occurs on an open subset of  $(M,g)$\,.\\ 
\indent
It is easy to see that a one-dim\-en\-sion\-al homothetic foliation has integrable orthogonal 
complement, at least outside the points where the foliation is Riemannian. Thus, if $\V$ 
is not Riemannian on $M$\,,   
we can find an open subset $U$ of $M$ on which $\V$ is not Riemannian and $\H$ is  integrable. Then $\O=0$ on $U$\,, and from  \eqref{e:riccixv} it follows that $X(V(\s))=(n-2)X(\s)V(\s)$ 
for any horizontal vector field $X$ where $\l=e^{\s}$ is the dilation of a harmonic morphism 
produced by $\V|_U$\,. By \eqref{e:Vinfautom} this is equivalent to 
\begin{equation} \label{e:xv} 
V(X((\s))=(n-2)X(\s)V(\s) 
\end{equation} 
for any basic vector field $X\in\G(\H)$\,.\\ 
\indent
Also, from \eqref{e:riccixy}\,, it follows that $X(\s)Y(\s)$ is a basic function for any orthogonal 
horizontal vectors $X$ and $Y$\,. Hence, by \eqref{e:xv}, we have that 
\begin{equation} \label{e:xy} 
0=V(X(\s)Y(\s))=2(n-2)X(\s)Y(\s)V(\s) 
\end{equation}     
for any orthogonal horizontal vectors $X$ and $Y$\,. Since $V(\s)\neq0$\,, it follows easily from   \eqref{e:xy} that, at each point $x\in U$ we have that $\H(\grad\s)_x=0$ (apply, 
for example, \cite[Lemma 3.2.3]{Pan-thesis}). Hence (see \cite{Pan-thesis}), case (ii) 
occurs on $U$ and the proof follows from Proposition \ref{prop:analytic}\,. 
\end{proof} 

\indent
Finally, we shall also need the following result from \cite[Proposition 3.4.8]{Pan-thesis}. 

\begin{prop} \label{prop:lateron} 
Let $(M,g)$ be an Einstein manifold, and let $\V$ be a one-dim\-en\-sion\-al foliation 
of codimension not equal to two which 
produces harmonic morphisms on $(M,g)$\,.\\
\indent
Then the following assertions are equivalent:\\
\indent
{\rm (i)} $\V$ has basic mean curvature form;\\
\indent
{\rm (ii)} $\V$ is a homothetic foliation. \qed
\end{prop}

\begin{rem} \label{rem:lateron} 
Let $\phi:(M,g)\to(N,h)$ be a submersive harmonic morphism. Let $\V={\rm ker}\,\phi_*$ 
and let $\l$ be the dilation of $\phi$\,. Then, from \eqref{e:BairdEells}\,, it follows that $\trace(\Bv)^{\flat}$ is basic if and only if $X(\log\l)$ is a basic function for any 
basic vector field $X\in\G(\H)$\,. 
\end{rem}

\indent
We now state the main theorem. 
\begin{thm} \label{thm:dmain} 
Let $(M,g)$ be an Einstein manifold of dimension at least $5$\,, and let $\V$ be a  
one-dim\-en\-sion\-al foliation which produces harmonic morphisms on $(M,g)$\,.\\ 
\indent
Then \emph{either},\\ 
\indent
$(\rm{i})$ $\V$ is a Riemannian foliation locally generated by Killing vector fields, \emph{or}\\ 
\indent
$(\rm{ii})$ $\V$ is a homothetic foliation by geodesics orthogonal to an umbilical foliation by hypersurfaces. 
\end{thm} 

\indent
Since by \cite{Bai} any harmonic morphism with one-dim\-en\-sion\-al fibres from a Riemannian  manifold 
of dimension at least 5 is submersive, from Theorem \ref{thm:dmain} we obtain the following.  

\begin{cor} \label{cor:dmain} 
Let $(M,g)$ be an Einstein manifold of dimension at least $5$\,, and let 
$\phi:(M,g)\to(N,h)$ be a harmonic morphism with one-dim\-en\-sion\-al fibres.\\ 
\indent
Then \emph{either},\\ 
\indent
$(\rm{i})$ the components of the fibres of $\phi$ form a Riemannian foliation locally tangent 
to nowhere zero Killing vector fields, \emph{or}\\ 
\indent
$(\rm{ii})$ $\phi$ is a horizontally homothetic submersion with geodesic fibres orthogonal to an 
umbilical foliation by hypersurfaces. 
\end{cor} 

\begin{rem} 
1) Note that, in case (i), $\phi$ is not, in general, a Riemannian submersion. In fact its dilation 
is proportional to $|V|^{1/(n-2)}$ where $V$ is a Killing vector field tangent to the fibres 
\cite{Bry} (see also \cite{BaiWoo2}\,,\,\cite{Pan}\,,\,\cite{Pan-thesis}\,).\\ 
\indent
2) Since, in case (ii), the hypersurfaces are the level hypersurfaces of $\l$ which is an isoparametric 
function, the fibres of $\phi$ are orthogonal to an isoparametric function \cite{Bai-book} 
(see \cite{BaiWoo2}\,). 
\end{rem} 

\begin{proof}[Proof of Theorem \ref{thm:dmain}]  
By Proposition \ref{prop:analytic}, if the horizontal distribution $\H$ is integrable on an open 
subset of $M$\,, then $\H$ is integrable on $M$\,. Thus, by Proposition \ref{prop:horintorhomoth}, 
it is sufficient to prove the case when $\H$ is nowhere integrable. Also, writing 
$n+1=\dim\,M$\,, we can suppose that 
the leaves of $\V$ are the fibres of a harmonic morphism 
$\phi:(M^{n+1},g)\to(N^n,\overline{h})$\,, 
where $\dim\,N=n$ with $n\geq3$\,, and from now on we shall use the notations of 
Lemma \ref{lem:Ricci}\,.\\ 
\indent
Let $\O^2\in\G({\rm End}(\H))$ be the field of self-adjoint negative semi-definite 
endomorphisms of $(\H,h|_{\H})$ defined by $h(\O^2(X),Y)=-h(i_X\O,i_Y\O)$ for horizontal 
$X$ and $Y$\,.\\ 
\indent
By Lemma \ref{lem:consdiag}\,, $\O^2$ can be consistently diagonalized on a dense open 
subset of $M$\,; let $x_0\in M$ be a point of this subset. 
There is an open neighbourhood $U$ of $x_0$ and an orthonormal frame $\{X_1,\ldots,X_n\}$    for $(\H,h|_{\H})$ over $U$ such that $\O^2(X_i)=-\mu_i^2\,X_i$ for some continuous 
functions $\mu_i:U\to[0,\infty)$ with $\mu_i^2$ smooth. Because $\O$ and $h|_{\H}$ are basic  
we also have that $\O^2$ is basic; hence the $\mu_i$ are basic as well. We can thus suppose 
that the $X_i$ are basic.\\ 
\indent
{}From \eqref{e:riccixy} we have 
\begin{equation} \label{e:1} 
\RicN(\phi_*X_i,\phi_*X_j)=(n-1)(n-2)X_i(\s)X_j(\s)\hspace{0.3cm}(i,j=1,\ldots,n\,,\:i\neq j).
\end{equation} 
\indent
We have the following alternative. \emph{Either}\\ 
(1) there exists $x\in U$ and distinct $j_1\,,\,j_2\,,\,j_3$ such that   
$X_{j_k}(\s)_x\neq0$\: $(k=1,2,3)$\,, \emph{or}\\ 
(2) for any $x\in U$ there are at most two distinct $j_1\,,\,j_2$ such that $X_{j_k}(\s)_x\neq0$\:   $(k=1,2)$\,.\\ 
\indent
Suppose that (1) holds. By \eqref{e:1} we have that $X_i(\s)X_j(\s)$ is basic for any 
$i\neq j$\,. Hence $X_{j_1}(\s)^2\,X_{j_2}(\s)^2\,X_{j_3}(\s)^2$ is basic, and, because 
$X_{j_k}(\s)\neq0$ on some open subset of $U$\,, we have that $X_{j_k}(\s)$ is basic\:   
$(k=1,2,3)$\,. Thus, if (1) holds, $X_i(\s)$ is basic for 
all $i=1,\ldots,n$\,, on some open subset of $U$\,. Then, by Proposition \ref{prop:analytic} 
and Proposition \ref{prop:lateron} (see also Remark \ref{rem:lateron}\,),  
$\V$ is homothetic on $M$ and the proof follows from Proposition \ref{prop:horintorhomoth}\,.\\ 
\indent
Suppose that (2) holds. If $X_j(\s)=0$ for all $j=1,\ldots,n$\,, then, by 
\eqref{e:BairdEells}\,, $\V$ has geodesic leaves and the proof of the theorem follows from 
Remark \ref{rem:homfoln}\,(3) and Proposition \ref{prop:horintorhomoth}\,.  
Therefore we can suppose that, after renumbering if necessary, we have $X_1(\s)_x\neq0$  
at some point $x\in U$\,. Then this holds on some open subset of $U$\,. Then, either 
$X_j(\s)=0$ for $j=2,\ldots,n$\,, on some open subset of $U$\,, or there 
exists a point $x\in U$ such that, after renumbering if necessary, $X_1(\s)_x\neq0$ 
and $X_2(\s)_x\neq0$\,. In the latter case, because (2) holds, we must have that 
$X_j(\s)=0$ ($j=3,\ldots,n$) on some open subset of $U$\,. It follows that there exists 
an open subset $U_1$ of $U$ such that $X_j(\s)=0$\,, ($j\geq3$). 
{}From now on we shall work on $U_1$\,.\\ 
\indent
By \eqref{e:riccixy} we have 
\begin{equation} \label{e:2} 
\begin{split}
c^M\,e^{-2\s}=\RicN(\phi&_*X_i,\phi_*X_i)\,-\,\tfrac12\,e^{(2n-2)\s}\,\mu_i^2\\  
&-e^{-2\s}\,\D^M\s\,-\,(n-1)(n-2)X_i(\s)^2\hspace{0.4cm}(i=1,\ldots,n)\,. 
\end{split}   
\end{equation} 
\indent
{}From \eqref{e:2} we get 
\begin{equation} \label{e:3} 
\begin{split} 
\RicN(\phi_*X_i,&\phi_*X_i)\,-\RicN(\phi_*X_j,\phi_*X_j)  
-\tfrac12\,e^{(2n-2)\s}\bigl(\mu_i^2-\mu_j^2\bigr)\\
&-\,(n-1)(n-2)\bigl(X_i(\s)^2-X_j(\s)^2\bigr)=0\hspace{0.4cm}(i,j=1,\ldots,n)\,. 
\end{split} 
\end{equation} 
\indent
{}From \eqref{e:3} it follows that $\frac12\,e^{(2n-2)\s}\bigl(\mu_i^2-\mu_j^2\bigr)$ 
is basic for $i,j\geq3$\,. Thus, if $\mu_i\neq\mu_j$ at some point for some $i,j\geq3$\,, 
$i\neq j$\,,  then $e^{\s}$ is basic and so $\V$ is Riemannian on some open subset of $M$\,;  
hence, by Proposition \ref{prop:analytic}\,, $\V$ is Riemannian on $(M,g)$\,. It 
remains to consider the case $\mu_3=\ldots=\mu_n=\mu$ for some function $\mu$\,.\\ 
\indent
Now, either, $\mu_1=\mu_2$ on some open subset, or, we have $\mu_1\neq\mu_2$
on a dense open subset. In the former case, by \eqref{e:3}, we have that 
$X_1(\s)^2-X_2(\s)^2$ is basic on some open subset. But, by \eqref{e:1}, 
$X_1(\s)X_2(\s)$ is also basic, and hence $X_1(\s)$\,, $X_2(\s)$ 
are basic on some open subset. Since $X_j(\s)=0$ for $j\geq3$\,, $\V$ has basic mean curvature form. 
Then, by Proposition \ref{prop:lateron}\,, 
$\V$ is homothetic on some open subset and hence, by Proposition \ref{prop:analytic}\,, 
$\V$ is homothetic on $(M,g)$\,; the proof of the theorem follows from 
Proposition \ref{prop:horintorhomoth}\,.\\ 
\indent
It remains to consider the case when $\mu_1\neq\mu_2$\,. Because $\O$ is skew-symmetric, at each 
point $x$\,, for 
any $i\in\{1,\ldots,n\}$ with $\mu_i(x)\neq0$ there exists $j\in\{1,\ldots,n\}$\,, $j\neq i$\,, 
such that $\mu_i(x)=\mu_j(x)$\,. Hence at each point $x$ we have that \emph{either}  
$\mu_1(x)=\mu(x)$ and $\mu_2(x)\neq\mu(x)$ \emph{or} $\mu_1(x)\neq\mu(x)$ 
and $\mu_2(x)=\mu(x)$\,. Suppose that $\mu_1(x)\neq\mu(x)$\,; then this holds  
at all points of an open subset, and on that subset we must have $\mu_2=\mu$\,. Moreoever, 
because $\O$ is skew-symmetric, 
we must have $\mu_1=0$ and so $\mu$ is not identically zero; in particular $n-1$ is even, 
i.e.,\ $n=2k+1$ for some integer $k\geq1$\,.\\ 
\indent
{}From \eqref{e:3} we get 
\begin{equation*} 
\RicN(\phi_*X_2,\phi_*X_2)-\RicN(\phi_*X_3,\phi_*X_3)\,=\,(n-1)(n-2)X_2(\s)^2\;; 
\end{equation*} 
hence, $X_2(\s)$ is basic. Thus, if $X_2(\s)\neq0$\,, since $X_1(\s)X_2(\s)$ is basic,   
we deduce that $X_1(\s)$ is also basic and the proof follows as before. There remains 
the case $X_2(\s)=0$ which we now consider. 
Summing-up the previous discussion, we have that $n=2k+1$\,, $k\geq1$\,, and we are now on an 
open subset on which we have the following: 
\begin{align*} 
&\mu_1=0\,,\:\mu_2=\ldots=\mu_n=\mu\,,\\ 
&X_2(\s)=\ldots=X_n(\s)=0\,,\\ 
&\mu\:{\rm and}\:X_1(\s)\:{\rm are\:not\:identically\:zero}.  
\end{align*}
Moreover, we can assume that $\mu$ and $X_1(\s)$ are \emph{nowhere} zero. Furthermore,  because 
$\O^2(X_1)=-\mu_1^2$ we have that $|i_{X_1}\!\O|_{\,h}^{\,2}=\mu_1^2=0$\,. Hence $i_{X_1}\!\O=0$\,, equivalently $i_{\grad_h\!\s}\O=0$\,.\\ 
\indent
{}From this and \eqref{e:riccixv} it follows that we have for $i=1,\ldots,n$\,, 
\begin{equation} \label{e:4}  
0=\tfrac12\,e^{(2n-2)\s}\,( ^h\dif^*\!\O)(X_i)+(n-1)X_i(V(\s))
-(n-1)(n-2)X_i(\s)V(\s)\;. 
\end{equation} 
\indent
Next, we compute $( ^h\dif^*\!\O)(X_1)$\,. 
\begin{equation*} 
\begin{split} 
( ^h\dif^*\!\O)(X_1)&=-\overset{2k+1}{\underset{j=1}{\sum}}
(\nabla_{X_j}\O)(X_j,X_1)\\ 
&=-\overset{2k+1}{\underset{j=1}{\sum}}\bigl\{X_j(\O(X_j,X_1)-\O(\nabla_{X_j}X_j,X_1)-\O(X_j,
\nabla_{X_j}X_1)\,\bigr\}\\ 
&=\overset{2k+1}{\underset{j=1}{\sum}}\O(X_j,\nabla_{X_j}X_1)\;=\; 
\overset{2k+1}{\underset{j=2}{\sum}}
\O(X_j,\nabla_{X_j}X_1)\\ 
&=\overset{k}{\underset{j=1}{\sum}}\bigl\{\,\O(X_{2j},\nabla_{X_{2j}}X_1)+
\O(X_{2j+1},\nabla_{X_{2j+1}},X_1)\,\bigr\}\;. 
\end{split} 
\end{equation*}  
We can choose a basic orthonormal local frame $\{X_1,X_2,\ldots,X_{2k+1}\}$ such that 
\begin{gather} \label{e:O}  
\left(\O_{ij}\right)=\begin{pmatrix} 
\,0 & 0 & 0 & \hdots & 0 & 0 \\ 
\,0 & 0 &-\mu& \hdots & 0 & 0 \\ 
\,0 &\mu& 0 &  \hdots & 0 & 0 \\ 
\hdotsfor{6}  \\ 
\,0 & 0 & 0 & \hdots & 0 & -\mu \\ 
\,0 & 0 & 0 & \hdots  & \mu & 0 
\end{pmatrix}\;.  
\end{gather}
Then, from the above calculation we have 
\begin{equation*} 
\begin{split} 
( ^h\dif^*\!\O)(X_1)&=\sum_{j=1}^{k}
\bigl\{h(\nabla_{X_{2j}}X_1,X_{2j+1})\,\O_{2j,2j+1}+ 
h(\nabla_{X_{2j+1}}X_1,X_{2j})\,\O_{2j+1,2j}\bigr\}\\ 
&=\sum_{j=1}^{k}\bigl\{-\mu\,h(\nabla_{X_{2j}}X_1,X_{2j+1})+\mu\,h(\nabla_{X_{2j+1}}X_1,X_{2j})\,
\bigr\}\\ 
&=\mu\sum_{j=1}^{k}\bigl\{h(X_1,\nabla_{X_{2j}}X_{2j+1})-h(X_1,\nabla_{X_{2j+1}}X_{2j})\,\bigr\}\\ 
&=\mu\sum_{j=1}^{k}h(X_1,[X_{2j},X_{2j+1}])\;.  
\end{split} 
\end{equation*} 
Recall that $X_j(\s)=0$ for all $j\geq2$\,; hence 
\begin{equation*} 
\begin{split} 
[X_{2j},X_{2j+1}](\s)=0\;&\iff\;-\V[X_{2j},X_{2j+1}](\s)=\H[X_{2j},X_{2j+1}](\s)\\ 
&\iff\;\O(X_{2j},X_{2j+1})\,V(\s)=h\bigl([X_{2j},X_{2j+1}],\H(\grad_h\s)\bigr)\\ 
&\iff\;-\mu\,V(\s)=h\bigl([X_{2j},X_{2j+1}],X_1\bigr)\,X_1(\s)\\ 
&\iff\;h\bigl([X_{2j},X_{2j+1}],X_1\bigr)=-\,\mu\,\frac{V(\s)}{X_1(\s)}\;. 
\end{split} 
\end{equation*}
It follows from the last equation that 
\begin{equation} \label{e:6} 
( ^h\dif^*\!\O)(X_1)=-\,k\,\mu^2\,\frac{V(\s)}{X_1(\s)}\;. 
\end{equation} 
\indent
{}From \eqref{e:4} and \eqref{e:6} we get 
\begin{equation*} 
0=-\tfrac12\,k\,\mu^2\,e^{4k\s}\,\frac{V(\s)}{X_1(\s)}+2k\,X_1(V(\s))-2k(2k-1)X_1(\s)V(\s) 
\end{equation*} 
which is equivalent to 
\begin{equation} \label{e:7} 
\mu^2\,e^{4k\s}\,V(\s)=4X_1(\s)X_1(V(\s))-4(2k-1)X_1(\s)^2\,V(\s)\;. 
\end{equation} 
\indent
{}From \eqref{e:3} with $i=1$\,, $j=2$\,, we get that 
$$\tfrac12\,e^{4k\s}\,\mu^2-2k(2k-1)X_1(\s)^2 \quad{\rm is\:basic}$$ 
and hence on differentiating this with respect to $V$\,, by \eqref{e:Vinfautom}\,, we obtain  
$$2k\,e^{4k\s}\,V(\s)\,\mu^2-4k(2k-1)X_1(\s)\,X_1(V(\s))=0$$ 
which is equivalent to 
\begin{equation} \label{e:8} 
\mu^2\,e^{4k\s}\,V(\s)=2(2k-1)X_1(\s)X_1(V(\s))\;.   
\end{equation} 
\indent
{}From \eqref{e:7} and \eqref{e:8} we get that 
$$4\,X_1(\s)X_1(V(\s))-4(2k-1)X_1(\s)^2\,V(\s)=2(2k-1)X_1(\s)X_1(V(\s))$$ 
which, because $X_1(\s)$ is nowhere zero, is equivalent to 
\begin{equation} \label{e:9} 
X_1(V(\s))=-\frac{2(2k-1)}{2k-3}\,X_1(\s)V(\s)\;. 
\end{equation} 
\indent
{}From \eqref{e:8} and \eqref{e:9} it follows that 
\begin{equation*}  
\mu^2\,e^{4k\s}\,V(\s)=-\frac{4(2k-1)^2}{2k-3}\,X_1(\s)^2\,V(\s) 
\end{equation*} 
which, if $\V$ is not Riemannian (equivalently, $V(\s)\neq0$), implies that 
\begin{equation} \label{e:10} 
\mu^2\,e^{4k\s}=-\frac{4(2k-1)^2}{2k-3}\,X_1(\s)^2\;.  
\end{equation} 
This is impossible if $k\geq2$\,, since $X_1(\s)\neq0$\,, $\mu\neq0$\,. The proof of the 
theorem is complete. 
\end{proof} 

\begin{rem} 
The same proof as above applies for the case $\dim M=4$ up to \eqref{e:10}\,.  
However, from \eqref{e:10} and $i_{X_1}\O=0$\,, we now have 
\begin{equation} \label{e:monopole} 
\dif\!^{\H}\!\left(\l^{-2}\right)=\ast_{\H}\,\O 
\end{equation} 
where $\dif^{\H}$ is the differential composed with the horizontal projection and $\ast_{\H}$ is 
the Hodge star-operator on $(\H,h|_{\H})$ with respect to some orientation  
of $\H$\,.\\ 
\indent
{}From the proof of Theorem \ref{thm:dmain} it follows that if $(M,g)$\,, $\dim M=4$\,, 
is an Einstein manifold and $\phi:(M,g)\to(N,h)$ is a submersive 
harmonic morphism with one-dim\-en\-sion\-al fibres which is not of type 1 or of type 2  
(i.e.\ $\V={\rm ker}\,\phi_*$ is neither Riemannian nor geodesic and with integrable 
horizontal distribution), then the `monopole equation' \eqref{e:monopole} must hold and  
we obtain \cite[Theorem 3.4.4]{Pan-thesis} (see also \cite{Pan-4to3}\,)\,. 
\end{rem}

\newpage 

\section{Harmonic morphisms with one-dim\-en\-sion\-al fibres between Einstein manifolds} 

In this section, we refine Theorem \ref{thm:dmain} in the case that the codomain is 
also Einstein. 

\begin{thm} \label{thm:dbetween} 
Let $\phi:(M,g)\to(N,h)$ be a surjective harmonic morphism with one-dim\-en\-sion\-al connected fibres  between Einstein manifolds, where $\dim M\geq5$\,. Then, up to homotheties, one of the following assertions 
holds:\\ 
\indent
{\rm (i)} $\phi$ is a Riemannian submersion with geodesic fibres 
onto an almost K\"ahler manifold and $\O$ is the pull-back of its K\"ahler 
form (in particular, $\dim M$ is odd).\\ 
\indent
{\rm (ii)} $\phi$ is a horizontally homothetic submersion with geodesic fibres orthogonal to an 
umbilical foliation by hypersurfaces.  
\end{thm} 

\begin{rem} 
Note that, for any $x\in M$\,, there exists an open neighbourhood $U$ of it, such that 
$\phi|_U$ has connected fibres. Furthermore, if $M$ is compact (more generally, if $\H$ is 
an Ehresmann connection) any harmonic morphism $\phi:(M^{n+1},g)\to(N^n,h)$ ($n\geq4$) 
can be factorised into a harmonic morphism with connected fibres followed by a Riemannian 
covering \cite{PanWoo-toprestr}\,. 
\end{rem} 

\begin{proof}[Proof of Theorem \ref{thm:dbetween}]
Let $n=\dim N$\,. If $\O=0$ then from \eqref{e:riccixy} it follows that 
$X(\s)Y(\s)=0$ for any pair 
of orthogonal vectors $X$ and $Y$ where $\l=e^{\s}$ is the dilation of $\phi$\,. 
Hence $X(\s)=0$ for any horizontal vector $X$\,, which implies that (ii) holds.\\ 
\indent
{}From now on we shall assume that $\O$ is nowhere zero. In particular (ii) 
does not hold and hence, by Theorem \ref{thm:dmain}\,, the 
fibres of $\phi$ form a Riemannian foliation locally generated by Killing vector 
fields. Therefore $\l=e^{\s}$ is basic and thus both $\l$ and $\s$ are  pull-backs 
by $\phi$ of functions on $N$\,; we shall denote these two functions 
by the same letters $\l$ and $\s$\,.\\ 
\indent 
We have $\RicM=c^Mg$\,, $\RicN=c^Nh$ and we shall denote by the same 
letter $h$ both the metric on $N$ and the metric on $M$ with respect to 
which $\phi$ is a Riemannian submersion with geodesic fibres.\\  
\indent
To complete the proof we must show that (i) holds on an open subset of $M$\,.\\  
\indent
{}From \eqref{e:riccixy} we get 
\begin{equation} \label{e:1between} 
\tfrac12\,e^{(2n-2)\s}\,h(i_X\O,i_Y\O)=-(n-1)(n-2)X(\s)Y(\s) 
\end{equation} 
for any pair $X,\,Y$ of orthogonal horizontal vectors.\\ 
\indent
As in the proof of Theorem \ref{thm:dmain}\,, there is an open subset $U$ of $M$ on which 
we can consistently diagonalise $\O^2$ with respect to a basic orthonormal frame $\{X_1,\ldots,X_n\}$ of $(\H,h)$\,. Then, from \eqref{e:1between} it follows that 
$X_i(\s)X_j(\s)=0$ for any $i\neq j$\,. It follows that either $X_i(\s)=0$ for all $i=1,\ldots,n$ 
or, after renumbering if necessary, on some open subset we have $X_2(\s)=\ldots=X_n(\s)=0$ and $X_1(\s)$ is not identically zero.\\ 
\indent
{}From \eqref{e:3} we obtain $\mu_i=\mu_j$ for any $i,j=2,\ldots,n$\,. 
Hence $\mu_i=\mu$ $(i=2,\ldots,n)$ for some non-negative function $\mu$ with $\mu^2$ smooth. 
{}From \eqref{e:3} we also get 
\begin{equation} \label{e:2between} 
-\tfrac12\,e^{(2n-2)\s}(\mu_1^2-\mu_2^2)=(n-1)(n-2)X_1(\s)^2\;. 
\end{equation} 
\indent
\emph{Suppose that} $X_1(\s)=0$ (equivalently, $\phi$ has geodesic fibres) then $\s$ is 
constant and 
$\mu_i=\mu$ for $i=1,\ldots,n$\,. Moreover, from \eqref{e:riccixy} it follows 
that $\mu=$\,constant. Hence $n$ must be even and $\O$ is $\mu$ times the K\"ahler form of 
an almost 
Hermitian structure on $(N,h)$\,. Since $\dif\!\O=0$\,, this structure is almost K\"ahler 
giving assertion (i)\,.\\ 
\indent
\emph{Suppose instead that} $X_1(\s)\neq0$\,. We shall obtain a contradiction.  
We have $\mu_1\neq\mu_2$ and hence, because 
$\O\neq0$ is skew-symmetric, $\mu_1=0$ and $\mu\neq0$\,. Thus,  
$h(i_{X_1}\O,i_{X_1}\O)=\mu_1^2=0$ which implies $i_{X_1}\O=0$\,. Together with $X_i(\s)=0$ $(i=2,\ldots,n)$\,, this gives that $i_{\grad\s}\O=0$\,. Also, \eqref{e:2between}  becomes 
\begin{equation} \label{e:3between} 
\tfrac12\,e^{(2n-2)\s}\,\mu^2=(n-1)(n-2)X_1(\s)^2 
\end{equation} 
which can be written as  
\begin{equation} \label{e:e} 
(n-1)\,\mu^2=2(n-2)\,X_1(\l^{1-n}\bigr)^2\;.  
\end{equation} 
\indent
Since $\phi$ is a harmonic morphism, by the chain rule (see, for example,  
\cite[Section 3.2]{BaiWoo2}\,), we have that 
$$\Delta\!^M(f\circ\phi)=\l^2\,\Delta\!^N\!f\circ\phi$$ 
for any smooth function $f$ on $N$\,.\\ 
\indent 
{}From \eqref{e:2} we get that 
\begin{equation} \label{e:4between} 
c^M\,e^{-2\s}=c^N-\Delta\!^N\!\s-(n-1)(n-2)X_1(\s)^2\;. 
\end{equation} 
\indent
{}From \eqref{e:riccivv} we get that 
\begin{equation} \label{e:5between} 
c^M\,e^{-2\s}=(n-2)\,\Delta\!^N\!\s+\tfrac14\,e^{(2n-2)\s}\,|\O|_h^2\;. 
\end{equation} 
Because $\mu\neq0$\,, as before, we have that $n$ must be odd. Then we can 
choose the frame $\{X_1,\ldots,X_n\}$ such that, with respect to this frame, $\O$ is 
given by  \eqref{e:O}\,. In particular, $|\O|_h^2=(n-1)\mu^2$ 
and \eqref{e:5between} becomes 
\begin{equation} \label{e:6between} 
c^M\,e^{-2\s}=(n-2)\,\Delta\!^N\!\s+\tfrac14(n-1)e^{(2n-2)\s}\mu^2\;. 
\end{equation} 
\indent
The relations \eqref{e:3between},\,\eqref{e:4between},\,\eqref{e:6between} give 
\begin{equation} \label{e:7between} 
\tfrac12(n-1)(n-2)(n-3)X_1(\s)^2=-(n-1)e^{-2\s}\,c^M+(n-2)c^N\;. 
\end{equation} 
Because $n\geq4$\,, relations \eqref{e:3between},\,\eqref{e:7between} and 
the fact that 
$X_i(\s)=0$ \;$(i=2,\ldots,n)$ imply $X_i(\mu)=0$ \;$(i=2,\ldots,n)$\,.\\ 
\indent
Note that \eqref{e:7between} can be written as follows: 
\begin{equation} \label{e:e1} 
\frac{n-3}{2}\,X_1(\l)^2=-\frac{c^M}{n-2}+\frac{c^N}{n-1}\,\l^2\;. 
\end{equation} 
\indent
{}From now on $\O$\,, $\mu$ and $X_j$\,, $j=1,\ldots,n$\,, will be viewed as  objects on $N$\,.  {}From  \eqref{e:riccixv} and $i_{\grad\s}\O=0$ we get that $^h\!\dif^*\!\O=0$\,.\\ 
\indent
Let $\F$ be the foliation formed on $N$ by the level hypersurfaces of 
$\s$\,. Note that $X_2\,,\ldots,\,X_n$ are tangent to $\F$ and $X_1$ is normal to $\F$\,. Also, 
$\mu$ is constant along the leaves of $\F$ and, from 
\eqref{e:O}\,, it follows that $\O$ restricted to any leaf $L$ of $\F$ 
is $\mu$ times the K\"ahler form of an almost Hermitian structure $J$ on 
$(L,h|_L)$\,. But $\dif\!\O=0$ and the constancy of $\mu$ along $L$ imply that 
$J$ is actually an almost K\"ahler structure on $(L,h|_L)$\,. In particular, 
$\Div_{\F}(\mu^{-1}\O|_L)=0$ which is equivalent to $\Div_{\F}(\O|_L)=0$ 
where $\Div_{\F}\O=\sum_{j=2}^n((\F\nabla)_{X_j}\O)(X_j,\cdot)$ 
(as is usual, we denote by the same letters the distributions 
$\F$ and $\F^{\perp}$ and the orthogonal projections onto them). 
After a short calculation (similar to the one used in the proof of 
Theorem \ref{thm:dmain} to find $(^h\!\dif^*\!\O)(X_1)$\,), by using the fact that $\F$ is integrable, we get that $(\Div_{\F}\O)(X_1)=0$\,. Hence  $\Div_{\F}\O=0$\,.\\ 
\indent
We claim that $\Div_{\F}\O=0$ and $\Div\O=-\,^h\!\dif^*\!\O=0$ imply 
that $\F^{\perp}$ has geodesic fibres. Indeed, we shall see that 
\begin{equation} \label{e:divO} 
\Div\O=\Div_{\F}\O-\mu\trace(\!B^{\perp}\!)^{\flat}\circ J 
\end{equation}
where $B^{\perp}$ is the second fundamental form of $\F^{\perp}$\,. Hence $\trace(\!B^{\perp}\!)=0$, i.e., $\F^{\perp}$ is geodesic.\\ 
\indent
Next, we prove \eqref{e:divO}\,. We have 
\begin{equation} \label{e:divO1} 
\begin{split} 
\Div\O&=\sum_{j=1}^n(\nabla_{X_j}\O)(X_j,\cdot)\\ &=(\nabla_{X_1}\O)(X_1,\cdot)
+\sum_{j=2}^n(\nabla_{X_j}\O)(X_j,\cdot)\\    
&=(\nabla_{X_1}\O)(X_1,\cdot)
+\sum_{j=2}^n((\F\nabla)_{X_j}\O)(X_j,\cdot)  
+\sum_{j=2}^n((\F\!^{\perp}\nabla)_{X_j}\O)(X_j,\cdot)\\ 
&=(\nabla_{X_1}\O)(X_1,\cdot)+\Div_{\F}\O+   
\sum_{j=2}^n((\F\!^{\perp}\nabla)_{X_j}\O)(X_j,\cdot)\;.   
\end{split}
\end{equation} 
Obviously, $(\nabla_{X_1}\O)(X_1,X_1)=0$\,, and $i_{X_1}\O=0$ implies that 
$((\F\!^{\perp}\nabla)_{X_j}\O)(X_j,X_1)=0$ for any $j=2,\ldots,n$ 
(recall that $\{X_1\}$ is a local frame for $\F^{\perp}$).  
Also $i_{X_1}\O=0$ and $X_j(\mu)=0$ \;$(j=2,\ldots,n)$ imply that 
$((\F\!^{\perp}\nabla)_{X_j}\O)(X_j,X_k)=0$ for any $j,k=2,\ldots,n$\,.\\ 
\indent 
For $j=1,\ldots,(n-1)/2$\,, we have 
\begin{equation} \label{e:divO2} 
\begin{split} 
(\nabla_{X_1}\O)(X_1,X_{2j})&=-\O(\nabla_{X_1}X_1,X_{2j})\\ 
 &=-h(\nabla_{X_1}X_1,X_{2j+1})\O(X_{2j+1},X_{2j})\\  
 &=-\mu\,h(\trace(\!B^{\perp}\!),JX_{2j})\;. 
\end{split} 
\end{equation} 
Similarly, 
\begin{equation} \label{e:divO3} 
(\nabla_{X_1}\O)(X_1,X_{2j+1})=-\mu\,h(\trace(\!B^{\perp}\!),JX_{2j+1})\;. 
\end{equation} 
\indent
Relation \eqref{e:divO} now follows from \eqref{e:divO1}\,,\,
\eqref{e:divO2} and \eqref{e:divO3}\,.\\ 
\indent
{}From $i_{X_1}\O=0$ and $\dif\!\O=0$ it follows that 
$\Lie_{X_1}\!\O=\dif i_{X_1}\O+i_{X_1}\!\dif\!\O=0$\,. Hence $\O$ is basic for  $\F^{\perp}$. It follows also that $J$ is basic 
for $\F^{\perp}$\,, equivalently, $\F^{\perp}$ is defined by horizontally 
holomorphic submersions (see the Appendix for the definition of  \emph{horizontally 
holomorphic submersion}). Moreover, by taking on the codomain of such a 
submersion the metric induced by any leaf of $\F$, from Proposition \ref{prop:horhol}  
of the Appendix we get that $\F^{\perp}$ 
is defined by horizontally holomorphic harmonic submersions.\\ 
\indent
It is easy to see that $\O(JX,Y)=\mu g(X,Y)$ for any $X,\,Y\in\F$\,. Because both 
$\O$ and $J$ are basic for $\F^{\perp}$ we get that $\F^{\perp}$ is a conformal 
foliation with dilation $\mu^{1/2}$\,. Moreover, because $\mu$ is constant along 
the leaves of $\F$ we get that $\F^{\perp}$ is a homothetic foliation. Thus, $\F^{\perp}$ 
is locally defined by horizontally homothetic submersions 
$\psi:(N,h)\to(P,k)$ with geodesic fibres (equivalently, $\F^{\perp}$ corresponds, 
locally, to a warped product decomposition of $(N,h)$\,). Moreover, $(P,k)$ is 
an Einstein $(n-2)$\,-manifold (apply, for example, 
\cite[Proposition 3.2.1]{Pan-thesis}\,)\,. {}From \cite[(3.2.1)]{Pan-thesis} (or 
\cite[(9.109)]{Bes}\,) it follows that $\mu$ satisfies the following equation: 
\begin{equation} \label{e:e2} 
X_1(\mu^{-1/2})=\frac{c^P}{n-2}-\frac{c^N}{n-1}\,\mu^{-1} 
\end{equation} 
where the Ricci tensor of $(P,k)$ is given by $\RicP=c^Pk$\,.\\ 
\indent
After a straightforward elementary calculation we get that the equations 
\eqref{e:e}\,, \eqref{e:e1} and \eqref{e:e2} are incompatible and the proof follows. 
\end{proof}

\begin{rem} 
1) Note that, if $\dim N=3$\,, then from the above proof we get that $i_{\grad\l}\O=0$ and \eqref{e:e}\,; equivalently, $\l^{-2}$ and $\O$ are related by the $S^1$-monopole 
equation (see \cite{Pan-thesis}, \cite{Pan-4to3} for this case).\\ 
\indent
2) Both case (i) and case (ii) of Theorem \ref{thm:dbetween} are related to 
well-known constructions of Einstein metrics (see \cite[Chapter 9]{Bes}\,).  
\end{rem}

\appendix

\section{A remark on horizontally holomorphic submersions} 

\indent
Let $\V$ be a foliation on the Riemannian manifold $(M,g)$\,. Assume that on the horizontal 
distribution $\H=\V^{\perp}$ there exists an almost Hermitian structure $J$ which is a 
basic tensor field for $\V$\,. Then, $\V$ could be locally defined by submersions 
$\phi:(U,g|_U)\to(N,\check{J})$ onto almost complex manifolds such that, at each point 
$x\in U$\,, 
the differential $\phi_{*,x}|_{\H_x}:(\H_x,J_x)\to(T_{\phi(x)}N,\check{J}_{\phi(x)})$ is complex 
linear. We shall call such submersions \emph{horizontally holomorphic} (see \cite{LouMo} 
for alternative terminology). Obviously, any holomorphic submersion from an almost Hermitian 
manifold to an almost complex manifold is horizontally holomorphic, but there exist simple examples of horizontally holomorphic submersions 
(e.g., orthogonal projections, the Hopf fibrations) which are not holomorphic  maps.\\ 
\indent 
Let $(M^m,g)$ be a Riemannian manifold endowed with a pair of complementary orthogonal 
distributions $\H$ and $\V$ of dimension $n$ and $m-n$\,, respectively. As is well-known, for a $(p,q)$\,-tensor field
$T$ on $M$ the \emph{horizontal divergence} $\Div_{\H}T$ 
is the  $(p,q-1)$\,-tensor field given by 
$$\Div_{\H}T=\sum_{a=1}^{n}i_{X_a}((\H\nabla)_{X_a}T)$$ 
where $\{X_a\}_{a=1,\ldots,n}$ is a local orthonormal frame for $\H$\,.  

\begin{prop} \label{prop:horhol} 
Let $\phi:(M,g)\to(N,h,J)$ be a horizontally holomorphic submersion to an almost 
Hermitian manifold. Then 
\begin{equation*} 
J(\trace(\Bv))+J(\tau(\phi))+\trace_g\phi^*(\N\!J)-\Div_{\H}\!J=0\;. 
\end{equation*} 
where  $\V={\rm ker}\,\phi_*$\,, $\Bv$ is the second fundamental form of $\V$,  $\N$ is the Levi-Civita connection of $(N,h)$ and $\tau(\phi)$ is the 
tension field of $\phi$\,. 
\end{prop} 
\begin{proof} 
As is well-known, by identifying $\phi^*(TN)$ with $\H$ the $\phi^*(TN)$\,-valued one-form induced 
by $\phi_*$ becomes the `horizontal' projection $\H$\,.\\ 
\indent
Let $X,\,Y$ be horizontal vector fields. With $\nabla$ the connection on 
$\phi^*(TN)\otimes T^*M$ we can write
\begin{equation*} 
\begin{split} 
(\nabla\phi_*)(X,JY)&=\phi^*(\N)_X(\H JY)-\H(\M_X(JY))\\ 
&=\phi^*(\N)_X(JY)-(\H\!\M)_X(JY)\\ 
&=(\phi^*(\N)_XJ)(Y)+J(\phi^*(\N)_XY)-((\H\!\M)_XJ)(Y)-J((\H\!\M)_XY)\;. 
\end{split} 
\end{equation*} 
Hence 
\begin{equation} \label{e:a1} 
(\nabla\phi_*)(X,JY)=J((\nabla\phi_*)(X,Y))+\phi^*(\N\!J)_XY-((\H\!\M)J)_XY\;. 
\end{equation} 
\indent
Let $\{X_a\}$ be a local orthonormal frame for $(\H,g_{\H})$\,. {}From \eqref{e:a1} we get 
\begin{equation} \label{e:a2} 
J\bigl(\,\sum_a(\nabla\phi_*)(X_a,X_a)\bigr)+\trace_g\phi^*(\N\!J)-\Div_{\H}\!J=0\;. 
\end{equation} 
\indent
As is well-known (and easy to prove) $(\nabla\phi_*)(V,W)=-\Bv(V,W)$ for any vertical vectors 
$V,\,W$\,. The proof follows from this fact and \eqref{e:a2}\,.  
\end{proof} 

\begin{rem} 
1) {}From Proposition \ref{prop:horhol} we could obtain a result similar to the one of  \cite{BaiEel} which relates the condition for a horizontally conformal  
submersion to be harmonic (and hence, a harmonic morphism) and the property that its 
fibres be minimal (\,\cite{LouMo}\,).\\ 
\indent
2) A proof similar to the above could be obtained for the following formula of A.~Lichnerowicz 
\cite{Lich} for a holomorphic map $\phi:(M,g,J)\to(N,h,J)$ between almost Hermitian 
manifolds: 
\begin{equation} \label{e:Lich} 
J(\tau(\phi))+\trace_g\phi^*(\N\!J)-\Div\!J=0\;. 
\end{equation} 
\indent
As is well-known from \eqref{e:Lich} it follows easily that any holomorphic map from 
a cosymplectic manifold to a (1,2)\,-symplectic manifold is harmonic \cite{Lich} 
(cf.\ \cite{GudWoo}\,). 
\end{rem}

\end{document}